\documentclass[12pt]{article}

\usepackage{amssymb, amsmath, amsthm, graphicx}
\usepackage[left=1in,top=1in,right=1in]{geometry}

\date{}

\theoremstyle{plain}

\theoremstyle{definition}

\theoremstyle{remark}

    \newcommand{\eps}{{\varepsilon}}

\begin{document}

\title{A stability theorem on cube tessellations}
\markright{Stability of cube tessellations}
\author{Peter Frankl\thanks{R\'enyi Institute, H-1364 Budapest, POB 127, Hungary. Email: {\tt peter.frankl@gmail.com}.} \and J\'anos Pach\thanks{R\'enyi Institute, Budapest, Hungary and EPFL, Lausanne, Switzerland. Email: {\tt pach@cims.nyu.edu}. Research partially supported by Swiss National Science Foundation Grants 200020-162884 and 200021-165977.}}

\date{}

\maketitle

\begin{abstract}
It is shown that if a $d$-dimensional cube is decomposed into $n$ cubes, the side lengths of which belong to the interval $(1-\frac{1}{n^{1/d}+1},1]$, then $n$ is a perfect $d$-th power and all cubes are of the same size. This result is essentially tight.
\end{abstract}

\section{Introduction}

It was proved by Dehn~\cite{De03} that, for $d\ge 2$, in any decomposition (tessellation, tiling) of the $d$-dimensional unit cube into finitely many smaller cubes, the side length of every participating cube must be rational. Fine and Niven~\cite{FN46} and, independently, Hadwiger raised the problem of characterizing, for a fixed $d\ge 2$, the set $N_d$ of all integers $n$ such that the $d$-dimensional unit cube can be decomposed into $n$ smaller cubes. Obviously, $m^d\in N_d$ for every positive integer $m$. Hadwiger observed that the intervals $(1,2^d)$ and $(2^d,2^d+2^{d-1})$ do not belong to $N_d$. On the other hand, for any $d$ there is a threshold $n_0(d)$ such that every integer $n\ge n_0(d)$ belongs to $N_d$; see \cite{P72}, \cite{M74}, \cite{E74}, \cite{CFG91}. It is conjectured that $n_0(d)\le c^d$ for a suitable constant $c$.

\smallskip

Amram Meir asked many years ago whether for any $d\ge 2, \eps>0$, and for every sufficiently large $n\ge n_0(d,\eps)$, there exists a decomposition of a $d$-dimensional cube into $n$ smaller cubes such that the ratio between the side lengths of any two cubes is at least $1-\eps$. This question was answered in the affirmative in~\cite{FMP17}. In particular, it was shown in~\cite{FMP17} that, for large $n$, a square can be decomposed into precisely $n$ smaller squares such that the ratio of their side lengths is at least  $1-O\left(\frac{1}{\sqrt{n}}\right)$.
\smallskip

The aim of this note is to show that the above bound is asymptotically tight. More precisely, we have the following stability result, which holds in every dimension $d\ge 2$.

\medskip

\noindent{\bf Theorem 1.} {\em Let $d, n\ge 2$ be positive integers. Suppose that a $d$-dimensional cube can be decomposed into precisely $n$ smaller cubes whose side lengths belong to the interval $(1-\frac{1}{n^{1/d}+1},1]$.

Then $n$ is a perfect $d$-th power, that is, $n=m^d$ for a positive integer $m$. Moreover, in this case the small cubes must be congruent.}
\medskip

\section{Proof of Theorem 1}

Consider a decomposition of the cube $[0,z]^d$ into $n$ smaller cubes of side lengths $s_i, 1\le i\le n$, where
$$1=s_1\ge s_2\ge \ldots \ge s_n>1-\frac{1}{n^{1/d}+1}.$$

By Dehn's theorem mentioned in the Introduction, we can assume that all $s_i$ and, hence, also $z$ are rational numbers. The total volume of the small cubes is $z^d$, so that we have
\begin{equation}\label{eq0}
z^d=\sum_{1\le i\le n}s_i^d\le ns_1^d =n.
\end{equation}
If equality holds here, then $s_1=\ldots=s_n=1$ and $n$ is a perfect $d$-th power, so we are done. Therefore, we can assume

\medskip

\noindent{\bf Claim 2.}  {\em $z < n^{1/d}.$}
\medskip

Fix a line $\ell$ parallel to the $x$-axis (say) that does not share a segment with the boundary of any small cube participating in the decomposition. (This holds, for example, if the other $d-1$ coordinates of the points of $\ell$ are all irrational.) Let $C_1, C_2, \ldots, C_m$ denote the small cubes crossed by $\ell$, listed from left to right, and let
$$0=x_0 < x_1 < x_2 <\ldots < x_m=z$$
be the $x$-coordinates of the points at which $\ell$ stabs the facets of these cubes. Using the assumption on the side lengths of the cubes, we have
\begin{equation}\label{formula1}
j\left(1-\frac{1}{n^{1/d}+1}\right)<x_j\le j,
\end{equation}
for every $j\; (1\le j\le m).$
\medskip

\noindent{\bf Claim 3.} {\em $m=\lceil z\rceil.$}
\medskip

\noindent{\bf Proof.} Since the side length of each cube $C_j$ is at most $1$, we clearly have $m\ge z$. It remains to show that $m<z+1$.

Suppose for contradiction that $m\ge z+1$. Applying (\ref{formula1}) with $j=m$, we obtain
$$z+\frac{n^{1/d}-z}{n^{1/d}+1}=(z+1)\left(1-\frac{1}{n^{1/d}+1}\right)
\le m\left(1-\frac{1}{n^{1/d}+1}\right)<x_m=z.$$
Comparing the left-hand side and the right-hand side, we get $n^{1/d}-z<0$, which contradicts Claim 2. $\Box$
\medskip

Claims 2 and 3 immediately imply that every line $\ell$ which is parallel to one of the coordinate axes and does not share a segment with the boundary of any small cube, intersects the same number, $m= \lceil z\rceil<n^{1/d}+1$, of small cubes. In particular, (\ref{formula1}) can be extended to
$$j-1<j\left(1-\frac{1}{n^{1/d}+1}\right)<x_j\le j,$$
for $1\le j\le m$. Thus, we can pick a small $\eps>0$ such that
\begin{equation}\label{formula2}
j-1+\eps \in (x_{j-1},x_j)
\end{equation}
holds for every $j\; (1\le j\le m)$.
\medskip

Given a small irrational number $\eps>0$, define a gridlike set $P_{\eps}$ of $m^d$ points in ${\mathbb R}^d$, as follows. Let
$$P_{\eps}=\{\eps, 1+\eps, 2+\eps,\ldots, m-1+\eps\}^d.$$
If $\eps$ is small enough, then all of these points lie in the interior of the cube $[0,z]^d$.
\medskip

\noindent{\bf Claim 4.} {\em  There exists $\eps>0$ such that every cube participating in the decomposition contains precisely one point in $P_{\eps}$.}
\medskip

\noindent{\bf Proof.} If $\eps$ is irrational, no element of $P_{\eps}$ lies on the boundary of any small cube. (This follows from the theorem of Dehn cited at the beginning of the Introduction.) The sidelength of every small cube is at most $1$, the minimum distance between two points in $P_{\eps}$, so that no cube can cover two elements of $P_{\eps}$.
\smallskip

We now finalize the choice of $\eps>0$. For every cube $C$ in the decomposition, pick a point $p=p(C)$ in the interior of $C$, all of whose coordinates are irrational. Let $\ell_1, \ell_2, \ldots, \ell_d$ denote the lines through $p$ parallel to the coordinate axes. None of them shares a segment with the boundary of any cube.
\smallskip

The line $\ell_1$ intersects precisely $m$ cubes. Suppose that $C$ is the $j$-th among them, and its projection to the first coordinate axis is the interval $[x_{j-1},x_j]$. If we choose $\eps>0$ small enough, then (\ref{formula2}) is satisfied for $\ell_1$. The same is true for the lines $\ell_2,\ldots,\ell_d$. Repeating the argument for every cube $C$, we can find an irrational $\eps>0$, which simultaneously satisfies all of the above conditions for all $C$. Then, for every $C$, there exist integers $j_k=j_k(C)$\; $(1\le j_k\le m,\, 1\le k\le d)$ such that the orthogonal projection of $C$ to the $k$-th coordinate axis contains $j_k-1+\eps$. Hence, we have
$$(j_1-1+\eps, j_2-1+\eps,\ldots, j_d-1+\eps)\in C,$$
showing that $C$ contains a point of $P_{\eps}$. $\Box$
\medskip

It follows from Claim 4 that $n$, the number of cubes participating in the decomposition, is equal to $|P_{\eps}|=m^d$. Thus, $n=m^d$ is a perfect $d$-th power.
\smallskip

Notice that the set $P_{\eps}$ can be covered by $m^{d-1}$ lines parallel to the first coordinate axis, and every small cube is stabbed by precisely one of these lines. The total sidelength of the cubes stabbed by each of these lines is equal to $z$. Therefore, the sum of the sidelengths of all small cubes satisfies $\sum_{i=1}^ns_i=\sum_{i=1}^{m^d}s_i=m^{d-1}z$, or, equivalently,
$$\frac{\sum_{i=1}^{m^d}s_i}{m^d}=\frac{z}{m}.$$

On the other hand, it follows from (\ref{eq0}) for $n=m^d$ that
$$\frac{\sum_{i=1}^{m^d}s_i^d}{m^d}=\left(\frac{z}{m}\right)^d.$$
For any positive numbers $s_i$, we have
$$\left(\frac{\sum_{i=1}^{m^d}s_i}{m^d}\right)^d\le
\frac{\sum_{i=1}^{m^d}s_i^d}{m^d},$$
with equality if and only if all $s_i$ are equal. In our setting equality holds, hence all small cubes must be of the same size.

This completes the proof of Theorem~1. \;\;\;\;\;\;\;\;\; $\Box$ $\Box$

\bigskip

Finally, we show that Theorem 1 is not far from being best possible. Consider the subdivision of the cube $[0,m]^d$ into $m^d$ unit cubes. Discard all of them that are not tangent to any of the coordinate hyperplanes. Fill out the resulting hole, $[1,m]^d$, by $m^d$ cubes of sidelength $1-\frac{1}{m}$. Altogether we have $$n=m^d-(m-1)^d+m^d<(m+1)^d=m^d+O(dm^{d-1})$$ 
cubes, where the inequality follows from the fact that the function $x^d$ is strictly convex.  The sidelengths of these cubes belong to the interval
$$[1-\frac{1}{m},1]=[1-\frac{1}{n^{1/d}(1+o(1))},1],$$
as $m$ tends to infinity. This interval is only slightly larger than the interval of ``permissible'' sidelengths in Theorem 1, but the number of small cubes participating in the tessellation is not a perfect $d$-th power.

\end{document}